\documentclass[journal]{IEEEtran}
\usepackage{xcolor,soul,framed} 
\colorlet{shadecolor}{yellow}
\usepackage[pdftex]{graphicx}
\graphicspath{{../pdf/}{../jpeg/}}
\DeclareGraphicsExtensions{.pdf,.jpeg,.png}
\usepackage[cmex10]{amsmath}
\usepackage{array}
\usepackage{mdwmath}
\usepackage{mdwtab}
\usepackage{eqparbox}

\usepackage{lineno}
\usepackage{amsmath}
\usepackage{graphicx}
\usepackage{float}
\usepackage{cite}
\usepackage{amssymb}
\usepackage{amsthm}
\usepackage{subfigure}
\usepackage{epsfig}
\usepackage{multirow}
\usepackage{pdfpages}
\usepackage{float}
\usepackage{nomencl}
\usepackage{subfigure}
\usepackage{epstopdf}
\usepackage{color}
\usepackage{url}
\usepackage{tabularx}
\usepackage{hhline}
\usepackage{rotating}
\usepackage{fancyhdr}
\usepackage[marginal]{footmisc}
\usepackage{microtype}
\usepackage{graphicx}
\usepackage{subfigure}
\usepackage{booktabs} 
\usepackage{amsmath}
\usepackage{multirow}
\usepackage{threeparttable}
\usepackage{stfloats}
\usepackage{makecell}

\usepackage[ruled,linesnumbered]{algorithm2e}

\usepackage{autobreak}

\usepackage{nomencl}

\usepackage{algpseudocode}
\algdef{SE}[DOWHILE]{Do}{doWhile}{\algorithmicdo}[1]{\algorithmicwhile\ #1}%

\allowdisplaybreaks[1] 

\usepackage{threeparttable}

\begin{document}
\title{Deep-Learning-Aided Voltage-Stability-Enhancing Stochastic Distribution Network Reconfiguration}
\author{Wanjun~Huang,~Changhong~Zhao,~\textit{Senior Member},~\textit{IEEE}
\thanks{This work was supported by the Hong Kong Research Grants Council through ECS Award No. 24210220. (Corresponding author: Changhong Zhao.)}
\thanks{W. Huang is with the School of Automation Science and Electrical Engineering, Beihang University, Beijing, China. Email: wjhuang1211@gmail.com}
\thanks{C. Zhao is with the Department of Information Engineering, the Chinese University of Hong Kong, New Territories, Hong Kong. Email: chzhao@ie.cuhk.edu.hk}
\thanks{This work has been submitted to the IEEE for possible publication. 
Copyright may be transferred without notice, after which this version may no longer be accessible.}
}
\maketitle

\newtheorem{thm}{Theorem} 
\newtheorem{definition}{Definition}
\newtheorem{lem}[thm]{Lemma}
\newtheorem{assumption}{Assumption}

\begin{abstract}
Power distribution networks are approaching their voltage stability boundaries due to the severe voltage violations and the inadequate reactive power reserves caused by the increasing renewable generations and dynamic loads. 
In the broad endeavor to resolve this concern, we focus on enhancing voltage stability through stochastic distribution network reconfiguration (SDNR), which optimizes the (radial) topology of a distribution network under uncertain generations and loads.
We propose a deep learning method to solve this computationally challenging problem.  
Specifically, we build a convolutional neural network model to predict the relevant voltage stability index from the SDNR decisions. 
Then we integrate this prediction model into successive branch reduction algorithms to reconfigure a radial network with optimized performance in terms of power loss reduction and voltage stability enhancement. 
Numerical results on two IEEE network models verify the significance of enhancing voltage stability through SDNR and the computational efficiency of the proposed method.
\end{abstract}


\section*{Nomenclature}
\addcontentsline{toc}{section}{Nomenclature}
\begin{IEEEdescription}[\IEEEusemathlabelsep\IEEEsetlabelwidth{$q_1,q_2$}]
\item[\textbf{Frequent acronyms:}]
\item[SDNR] Stochastic distribution network reconfiguration.
\item[RVSI] Root-mean-squared voltage-dip severity index (for short-term voltage stability). 
\item[SBR] Successive branch reduction.
\item[CNN] Convolutional neural network.
\\

\item[\textbf{Sets:}]
\item[$\mathcal{N}$] The set of buses in a distribution network, including substation buses $\mathcal{N}_s$ and non-substation buses $\mathcal{N}_d$.
	\item[$\mathcal{E}$] The set of branches $e=ij \in\mathcal{E}$.
	\item[$\mathcal{W}$] The set of scenarios $w\in\mathcal{W}$ for uncertain renewable generations and loads. 
	\item[$\mathcal{A}$] The set of feasible switch status vectors $\boldsymbol{\alpha}\in\mathcal{A}$ that lead to radial networks.
	\item [$\mathcal{P}$] The set of branches that form a loop.
	\item [$\mathcal{K}$] The set of candidate branches to open.
 \\

	\item[\textbf{Given quantities:}] 
	\item[$\boldsymbol{\pi}$] The probability distribution $\boldsymbol{\pi}=(\pi_w,\forall w\in\mathcal{W})$ of the uncertainty scenarios. 

\item[$\hat{p}_{ri}^{w},\hat{q}_{ri}^{w}$] The active and reactive renewable power generations at bus $i$ in scenario $w$. 
	
 \item[$\hat{p}_{di}^{w},\hat{q}_{di}^{w}$] The active and reactive power loads at bus $i$ in scenario $w$.

	\item[$g_{ij}, b_{ij}$] Series conductance and susceptance of branch $ij$.

 \item[$L$] The number of redundant branches in a network, which form the same number of chordless loops.
\\ 

\item[\textbf{Variable quantities:}]
	\item[$C_{l}^{w}$] The total active power loss in scenario $w$.
    \item[$I_{v}^{w}$] The voltage stability index in scenario $w$.
    \item [$\delta_{\text{min}}$] The smallest singular value of the power-flow Jacobian matrix (indexing steady-state voltage stability).

\item[$\alpha_{ij}$] The binary variable indicating the switch status on branch $ij$, collected in $\boldsymbol{\alpha}=(\alpha_{ij},\forall ij\in\mathcal{E})$.
    
	\item[$p_{i}^{w},q_{i}^{w}$] The active and reactive power injections at bus $i$ in scenario $w$.
	
	\item[$V_{i}^{w},\theta_{i}^{w}$] The voltage magnitude and phase angle at bus $i$ in scenario $w$.

 \item[$p_{ij}^w,q_{ij}^{w}$] The active and reactive power flows on branch $ij$ in scenario $w$.
 
	\item[$\boldsymbol{x}$] The continuous decision variables:
	\begin{eqnarray}
\boldsymbol{x}&=&\left(\boldsymbol{x}^{w},~\forall w\in\mathcal{W}\right)\nonumber\\
	&=& \left(p_{i}^{w}, q_{i}^{w}, \forall i\in \mathcal{N}_s;\right.\nonumber\\
	&& \left.~V_{i}^{w},\theta_{i}^{w},\forall i\in\mathcal{N};~\forall w\in\mathcal{W}\right).\nonumber
	\end{eqnarray}
	\item[$\widetilde{p}_{ij}$] The expected active power flow on branch $ij$.
	
	\item[$\widetilde{p}_{i}^{\mathcal{P}}$] The expected active power injection from bus $i$ into loop $\mathcal{P}$. 

	\item[$n_r$] The number of active power-injecting buses that divide loop $\mathcal{P}$ into sub-paths.

	
\end{IEEEdescription}

\section{Introduction}
Voltage stability is becoming a major issue in power distribution networks
with increasing shares of renewable generations and dynamic loads. 
A major cause of voltage instability is the lack of reactive power supplies, as the inverters that interface the renewable energy sources 
are demanding more reactive power compensation \cite{bajaj2020grid}. 
Meanwhile, the larger variations in distributed renewable generations and loads, as well as the increasingly severe contingencies from the upper/external grid \cite{federal2012arizona}, are pushing distribution networks closer to their operating limits, particularly voltage stability boundaries.   
Therefore, it is critical to enhance voltage stability in the operations of renewables-heavy distribution networks.

It has been recognized that the network topology plays a significant role in voltage stability \cite{simpson2016voltage,
huang2020small}.
The network topology can be optimized through distribution network reconfiguration (DNR), which configures the opening/closure of switches on the branches (power lines) to improve network performance. 
Typical DNR formulations mainly focused on power loss minimization and load balancing 
\cite{baran1989network}.
Other than that, voltage stability is an objective for DNR that should receive more attention than it currently does.

The DNR problem, even without considering voltage stability, is known to be a difficult mixed-integer nonlinear program due to the binary switching actions subject to the nonconvex power flow constraints and the requirement for a radial (tree) topology. 
Popular solution methods for DNR encompass mathematical programming, heuristics or meta-heuristics, and machine learning, each facing certain challenges in practice. 
For instance, mathematical programs that exploit various convex approximations \cite{Taylor2012convex,Jabr2021minimum} need to address non-trivial trade-offs between accuracy, optimality, and computational efficiency; heuristics such as the iterative branch exchange \cite{civanlar1988distribution} and switch opening and exchange \cite{zhan2020switch} may have sensitive performance to initialization and poor convergence rates; meta-heuristics such as the genetic algorithm \cite{jakus2020optimal} and particle swarm optimization \cite{pegado2019radial} may suffer heavy computations and inconsistent outcomes from different randomized runs. Compared to the methods above, a class of successive branch reduction (SBR) heuristics   \cite{peng2014feeder,huang2022improved} can reach a satisfactory balance between optimality and computational efficiency. 
     

Taking voltage stability into account will make the DNR problem even harder. 
Most prior efforts for voltage-concerned DNR aimed to flatten voltage profiles \cite{imran2014novel,nguyen2015distribution,wu2020distribution} or restrict voltage volatility \cite{song2020new}. 
Fewer studies tried to improve steady-state voltage stability \cite{nguyen2016novel,raut2020improved,tran2021reconfiguration} and small-disturbance voltage stability \cite{shukla2020stochastic} through DNR.
The more complicated indices for (large-disturbance) short-term voltage stability are generally calculated from real-world or simulated time-series voltage records rather than expressed explicitly in terms of the DNR and volt/var control decisions \cite{zhang2019hierarchical}.
This makes it difficult to integrate the voltage stability indices into the DNR formulations and to optimize them through the traditional mathematical programming, heuristic, or meta-heuristic approaches.   
This difficulty is partly addressed by the deep learning method we recently developed to predict short-term voltage stability from DNR decisions \cite{huang2021distribution,huang2021resilient}.


It is worth noting that the robust or stochastic version of DNR, rather than the simple deterministic version, need to be implemented to deal with the uncertain scenarios of renewable generations and loads.
In this paper, we focus on the stochastic DNR (SDNR), which optimizes the expectation across such scenarios and is thus generally less conservative and more economic than the robust version 
that concerns about the worst scenario.
 Grounded in our recent SBR heuristics for SDNR \cite{huang2022improved} and our deep learning method to predict voltage stability in deterministic DNR \cite{huang2021distribution}, we introduce the following new method and result to this research field:
\begin{itemize}
    \item A deep learning method is proposed to solve SDNR with voltage stability enhancement. First, a convolutional neural network (CNN) model is built to predict the voltage stability index concerned from SDNR decisions. 
    Then, the CNN model is integrated into the one-stage and two-stage SBR algorithms from \cite{huang2022improved} to reconfigure a radial network with optimized performance in terms of power loss reduction and voltage stability enhancement.  
   
   \item The proposed method is validated by numerical experiments on two IEEE distribution network models. Compared to a classic mixed-integer nonlinear program solver for SDNR without considering voltage stability, the proposed method improves the steady-state or short-term voltage stability at a minor cost of increased power loss. Moreover, it speeds up the solution process of SDNR by at least an order of magnitude.
\end{itemize}

In the following, Section \ref{sec:model} introduces our model and formulation of voltage-stability-enhancing SDNR. Section \ref{sec:proposed_method} elaborates on our deep-learning-aided method to solve the formulated SDNR problem. The numerical case studies are presented in Section \ref{Sec:CaseStudy}. Then Section \ref{sec:conclusion} concludes the paper.

\section{Model and Problem Formulation}\label{sec:model}

\subsection{SDNR with Voltage Stability Enhancement}
Similar to our settings in \cite{huang2022improved}, we consider a distribution network with a set of buses $\mathcal{N}$ and a set of branches $\mathcal{E}$. 
The bus set $\mathcal{N}$ is divided into substation buses $\mathcal{N}_{s}$ and non-substation buses $\mathcal{N}_{d}$.
The substation buses are connected to an upper-level transmission network, while the non-substation buses are connected to loads and/or renewable energy sources. 

Each branch in the set $\mathcal{E}$ is arbitrarily assigned a reference direction, say from bus $i$ to bus $j$, and is denoted as $ij \in\mathcal{E}$; the existence of $ij\in\mathcal{E}$ shall exclude the other direction $ji\in \mathcal{E}$, so that the branches are not double counted.
In case the quantities associated with both directions of a branch need to be used, we define the unordered relationship $i \sim j$ and $j\sim i$, which \textit{simultaneously} hold as long as there is a branch between buses $i$ and $j$ in \textit{either} direction.
Without loss of generality, we assume all the branches $ij \in\mathcal{E}$ are switchable. A binary variable $\alpha_{ij}$ indicates the open ($\alpha_{ij}=0$) and closed ($\alpha_{ij}=1$) status of the switch on branch $ij$.

The uncertain renewable generations and loads are modeled as random variables subject to a joint probability distribution $\boldsymbol{\pi}=(\pi_w,\forall w\in \mathcal{W})$ over a finite set of scenarios $\mathcal{W}$. 

With the settings above, we formulate an SDNR problem with voltage stability enhancement as follows.

\subsubsection{Objective} We aim to optimize the expected operational cost across all the scenarios $w\in\mathcal{W}$ in terms of power loss and voltage stability:
\begin{eqnarray}
\Gamma(\boldsymbol{\alpha},\boldsymbol{x}) &=&\sum_{w\in\mathcal{W}}\pi_w\left(k_{l}\frac{C_{l}^{w}}{C^{\text{max}}_{l}} + k_{v}\frac{I_{v}^{w}}{I^{\text{max}}_{v}}\right) 
    \label{eq:SDNR_obj}
\end{eqnarray}
where the decision variables $(\boldsymbol{\alpha},\boldsymbol{x})$ will become clear as we proceed. 
The total active power loss $C_{l}^{w}$ and voltage stability index $I_{v}^{w}$ are normalized by their upper bounds (set as constants by experience), $C^{\text{max}}_{l}$ and $I^{\text{max}}_{v}$, respectively. 
A desired balance can be reached between the two objectives, by tuning the factors $k_{l}$ and $k_{v}$.
Note that $k_{l}$, $C_{l}^{w}$, ${C^{\text{max}}_{l}}$, $I_{v}^{w}$, ${I^{\text{max}}_{v}}$ are always positive numbers, while $k_{v}$ can be positive or negative depending on whether a smaller or larger $I_{v}^{w}$
is more stable for the specific type of voltage stability concerned, as will be explained in Section~\ref{subsec:voltage_stability}. The formulation of $C_{l}^{w}$ will be introduced shortly with the bus power injections.

\subsubsection{Branch power flows} 
The active and reactive branch power flows are $\forall ij\in \mathcal{E}$, $\forall w\in\mathcal{W}$: 
\begin{subequations}\label{eq:branch_flow}
\begin{eqnarray}
p_{ij}^{w}=\alpha_{ij}\!\left[(V^{w}_{i})^{2}g_{ij}\!-\!V_{i}^{w}V_{j}^{w}(b_{ij}\text{sin}\theta_{ij}^{w}\!+\!g_{ij}\text{cos}\theta_{ij}^{w})\right] &&\label{eq:pij}\\
p_{ji}^{w}=\alpha_{ij}\!\left[(V^{w}_{j})^{2}g_{ji}\!-\!V_{j}^{w}V_{i}^{w}(b_{ji}\text{sin}\theta_{ji}^{w}\!+\!g_{ji}\text{cos}\theta_{ji}^{w})\right] && \label{eq:pji} \\
 q_{ij}^{w}=\alpha_{ij}\!\left[-(V^{w}_{i})^{2}b_{ij}\!-\!V_{i}^{w}V_{j}^{w}(g_{ij}\text{sin}\theta_{ij}^{w}\!-\!b_{ij}\text{cos}\theta_{ij}^{w})\right] &&\label{eq:qij} \\
 q_{ji}^{w}=\alpha_{ij}\!\left[-(V^{w}_{j})^{2}b_{ji}\!-\!V_{j}^{w}V_{i}^{w}(g_{ji}\text{sin}\theta_{ji}^{w}\!-\!b_{ji}\text{cos}\theta_{ji}^{w})\right] &&\label{eq:qji}
\end{eqnarray}
\end{subequations}
where $V_i^w$, $V_j^w$ are the voltage magnitudes at buses $i$, $j$ and  $\theta_{ij}^{w}:=\theta_{i}^{w}-\theta_{j}^{w}$ is the angle difference between them.  
Due to power loss, we cannot say $p_{ij}^w = -p_{ji}^w$, $q_{ij}^w = -q_{ji}^w$, i.e., the power sent onto a branch from one end generally does not equal the power received by the other end. The power at both ends are calculated in \eqref{eq:branch_flow}, because they will both appear in constraint \eqref{eq:power_balance} below that exploits the unordered relationship $i\sim j$.
The series admittance of a branch is $g_{ij} + \mathrm{j} b_{ij} = g_{ji} + \mathrm{j} b_{ji}$ regardless of its direction. For each $ij\in\mathcal{E}$, it is sufficient to define just one switch status $\alpha_{ij} \in\{0,1\}$ to determine all the four branch flows in \eqref{eq:pij}--\eqref{eq:qji}.

\subsubsection{Bus power injections} 
The power injection at each bus balances the total power flow sent to its adjacent buses:
\begin{eqnarray}
p_{i}^{w} = \sum_{j: i\sim j} p_{ij}^w, \quad q_{i}^{w} = \sum_{j: i\sim j} q_{ij}^w, ~\forall i\in\mathcal{N},~\forall w\in\mathcal{W}. \label{eq:power_balance}
\end{eqnarray} 
The power imported at the substation buses are restricted as:
\begin{eqnarray}
p_{i}^{\text{min}}\!\leq\! p_{i}^{w} \!\leq\! p_{i}^{\text{max}},~ q_{i}^{\text{min}}\!\leq\! q_{i}^{w} \!\leq \!q_{i}^{\text{max}},~ \forall i\in\mathcal{N}_s,~\forall w\in\mathcal{W}\label{eq:gen_limit} 
\end{eqnarray}
where $(\cdot)^{\text{min}}$ and $(\cdot)^{\text{max}}$ are the given minimum and maximum limits. The power injections at non-substation buses are:
\begin{eqnarray}
p_{i}^{w}=\hat{p}_{ri}^{w}-\hat{p}_{di}^{w},~ q_{i}^{w}=\hat{q}_{ri}^{w}-\hat{q}_{di}^{w},~\forall i\in\mathcal{N}_{d},~\forall w\in\mathcal{W}
\label{eq:bus_inj}
\end{eqnarray}
where $(\hat{p}_{ri}^{w},\hat{q}_{ri}^{w})$ are the active and reactive power generation of the aggregate renewable energy source, and $(\hat{p}_{di}^{w},\hat{q}_{di}^{w})$ are the active and reactive power consumption of the aggregate load, at bus $i$. As mentioned, they are uncertain quantities subject to probability distribution $\boldsymbol{\pi}$ over scenarios $w\in\mathcal{W}$. 

Summing up the active power injections at all the buses leads to the total active power loss in objective \eqref{eq:SDNR_obj}:
\begin{eqnarray}
C_{l}^{w}&:=&\sum_{i\in\mathcal{N}} p_i^w.\nonumber
\end{eqnarray}

\subsubsection{Safety limits} The bus voltages are limited as: 
\begin{eqnarray}
&&V_{i}^{\text{min}}\leq V_{i}^{w} \leq V_{i}^{\text{max}},~\forall i \in \mathcal{N},~\forall w\in\mathcal{W} \label{eq:V}
\end{eqnarray}
and the apparent power flows are limited as:
\begin{eqnarray}
&&(p_{ij}^{w})^2+(q_{ij}^{w})^2\leq (s_{ij}^{\text{max}})^2,~\forall i\sim j,~\forall w\in\mathcal{W}. \label{eq:sijmax}
\end{eqnarray}

\subsubsection{Network topology}
Denote the network graph under switch status $\boldsymbol{\alpha}=(\alpha_{ij},\forall ij\in \mathcal{E})$ as $\mathcal{G}(\boldsymbol{\alpha})$. The feasible set of $\boldsymbol{\alpha}$ is defined as \cite{peng2014feeder}:
\begin{eqnarray}
\mathcal{A}:=\{\boldsymbol{\alpha}~ |&\mathcal{G}(\boldsymbol{\alpha})~\textnormal{has no loop; and each bus in}~\mathcal{N}_d \nonumber\\
&\textnormal{is connected to a single bus in}~\mathcal{N}_s\}. \nonumber
\end{eqnarray}
A network $\mathcal{G}(\boldsymbol{\alpha})$ is \textit{radial} if it satisfies the condition in $\mathcal{A}$. 

To reduce the variable space, we treat $(p_{ij}^w,p_{ji}^w, q_{ij}^w, q_{ji}^w)$ in \eqref{eq:branch_flow} as expressions rather than variables. 
Now we are ready to introduce the SDNR with voltage stability enhancement:
\begin{eqnarray} \label{eq:SDNRVS}
	&&\textbf{SDNR-VS}:~\mathop{\text{minimize}}~\Gamma(\boldsymbol{\alpha},\boldsymbol{x})\nonumber\\
	&&\quad\quad\quad \text{over}\quad \boldsymbol{\alpha}:=(\alpha_{ij},~\forall ij\in \mathcal{E})\in \mathcal{A}, \nonumber\\
	&&\quad\quad\quad\quad~~~~~ \boldsymbol{x}:=(\boldsymbol{x}^{w},~\forall w\in\mathcal{W}) \nonumber \\
&& := \left(p_{i}^{w}, q_{i}^{w},\forall i\in \mathcal{N}_s;~V_{i}^{w},\theta_{i}^{w},\forall i\in\mathcal{N};~\forall w\in\mathcal{W}\right)\nonumber\\
&&\quad\textnormal{subject to~~ \eqref{eq:branch_flow}--\eqref{eq:sijmax}}. \nonumber
\end{eqnarray}

For each given and fixed switch status $\boldsymbol{\alpha}$, the \textbf{SDNR-VS} problem is specified as a stochastic optimal power flow (OPF) problem taking voltage stability into account. 
As part of the proposed procedure to solve \textbf{SDNR-VS} in Section \ref{sec:proposed_method}, we will deal with a reduced version of this stochastic OPF problem, which minimizes the expected total active power loss only, without considering voltage stability:
\begin{eqnarray}
	\textbf{SOPF}(\boldsymbol{\alpha}):&&\mathop{\text{minimize}}\limits_{\boldsymbol{x}}~\sum_{w\in\mathcal{W}}\pi_w C_{l}^{w}\nonumber\\
	&&\text{subject to}\quad \text{\eqref{eq:branch_flow}--\eqref{eq:sijmax}}.  \nonumber
\end{eqnarray}	


Classic SDNR problems are mixed-integer nonlinear programs solvable by off-the-shelf software such as Gurobi. 
Compared to them, a major hurdle in our formulation \textbf{SDNR-VS} is that the voltage stability index $I_{v}^{w}$ in objective \eqref{eq:SDNR_obj} may not have an explicit expression in terms of the decision variables. For instance, the large-disturbance short-term voltage stability index is generally calculated from real-world or simulated time-series voltage records. 
Therefore, it will be difficult for the existing mathematical program solvers to handle \textbf{SDNR-VS}. This motivates the deep-learning-aided method in this paper.

\subsection{Voltage Stability Evaluation}\label{subsec:voltage_stability}

The deep learning method to be proposed shortly can predict the voltage stability index $I_{v}^{w}$ from the SDNR decisions, for various types of voltage stability  defined in \cite{hatziargyriou2021definition}. In this paper, we just focus on the steady-state and short-term voltage stability as representative examples.

\subsubsection{Steady-state voltage stability} 
At a steady-state operating point (i.e., a power flow solution), the smallest singular value of the Jacobian matrix, denoted by $\delta_{\text{min}}$, indicates the margin of this operating point from a voltage collapse \cite{lim2016svd}. As $\delta_{\text{min}}$ approaches zero, even a small variation in active or reactive power injection will cause a large voltage excursion.  
Therefore, one can apply the singular value decomposition method to calculate $\delta_{\text{min}}$ as a stability index. 

\subsubsection{Short-term voltage stability} It focuses on the dynamic behavior of a network after a large disturbance such as a short-circuit fault. 
It is typically assessed via the time-series simulation of voltage dynamics \cite{zhang2019hierarchical,
huang2021resilient}, which often involves induction motor and ZIP load models \cite{zhang2020identification}.   
The short-term voltage stability problems generally encompass two issues: the voltage instability and the fault-induced delayed voltage recovery \cite{huang2021distribution}. 
The former can be directly identified from the time-series voltage records. 
The latter can be evaluated by the root-mean-squared voltage-dip severity index (RVSI), which should be smaller for a more stable network \cite{zhang2019hierarchical}. 

\section{Deep-Learning-Aided Solution Method}\label{sec:proposed_method}

\subsection{CNN-Based Prediction of Voltage Stability}\label{subsec:CNN}

As the basis of the proposed algorithms to solve the challenging \textbf{SDNR-VS} problem, we briefly introduce the convolutional neural network (CNN) model from our recent work \cite{huang2021distribution} to predict voltage stability. Through proper training, a CNN can learn the nonlinear mapping 
\begin{eqnarray}
F_{vs}:~\boldsymbol{G}(\boldsymbol{\alpha},\boldsymbol{x})\rightarrow I_{v} \nonumber
\end{eqnarray}
that maps a network state to its corresponding voltage stability index. 
Here $\boldsymbol{G}(\boldsymbol{\alpha},\boldsymbol{x})$ collects the network state as
\begin{eqnarray}
\boldsymbol{G}(\boldsymbol{\alpha},\boldsymbol{x}) := \left(\boldsymbol{\ell}_{ij}, ~\forall ij\in \mathcal{E}_{\boldsymbol{\alpha}}\right) \nonumber
\end{eqnarray} 
where $\mathcal{E}_{\boldsymbol{\alpha}} := \{ij\in \mathcal{E}~|~\alpha_{ij}=1\}$ is the set of closed branches determined by the switch status vector $\boldsymbol{\alpha}$, and vectors
\begin{eqnarray}
\boldsymbol{\ell}_{ij} := \left[i, j, g_{ij}, b_{ij}, p_{ij}, q_{ij}, p_{j}, q_{j}\right]^{\intercal}, ~\forall ij\in\mathcal{E}_{\boldsymbol{\alpha}} \nonumber
\end{eqnarray}
contain the parameters and variables associated with the closed branches.
The output $I_{v}$ can be any voltage stability index of the operator's interest, e.g., the smallest singular value of the power-flow Jacobian matrix for steady-state voltage stability or the RVSI for short-term voltage stability, as already reviewed in Section \ref{subsec:voltage_stability}. 

\subsection{Deep-Learning-Aided One-Stage SBR Algorithm}

The proposed algorithms to solve \textbf{SDNR-VS} extend the successive branch reduction (SBR) heuristics in our recent work \cite{huang2022improved} to incorporate the deep-learning-aided prediction of voltage stability. 
When introducing the proposed algorithms below, we shall just highlight their key differences and skip their common details and rationales with those in \cite{huang2022improved}. 

We assume the distribution network has a single substation bus in $\mathcal{N}_s$, since merging multiple substations into one would not change the applicability of the proposed algorithms. 
First, we present a one-stage SBR algorithm, Algorithm \ref{alg:OSSBR}, for a simple special network with a single redundant branch than radial. We then present a two-stage SBR algorithm, Algorithm \ref{alg:TSSBR}, that can be applied to a general network with multiple redundant branches. In its second stage, Algorithm \ref{alg:TSSBR} will iteratively call Algorithm \ref{alg:OSSBR}.  

\begin{algorithm}[t]
	\caption{Deep-Learning-Aided One-Stage SBR}\label{alg:OSSBR}
	\SetKwData{Left}{left}\SetKwData{This}{this}\SetKwData{Up}{up}\SetKwFunction{Union}{Union}\SetKwFunction{FindCompress}{FindCompress}\SetKwInOut{Input}{Input}\SetKwInOut{Output}{Output}
	Initialize switch status as $\boldsymbol{\alpha}_{\mathcal{E}}:=(\alpha_{ij}=1,\forall ij\in\mathcal{E})$;\\
	Solve $\textbf{SOPF}(\boldsymbol{\alpha}_{\mathcal{E}})$ to get optimal $\boldsymbol{x}_\mathcal{E}=(\boldsymbol{x}^w_\mathcal{E},\forall w\in\mathcal{W})$;\\
	At $\boldsymbol{x}_\mathcal{E}$, calculate $\widetilde{p}_{i}^{\mathcal{P}}$ by \eqref{eq:pi_equivalent} for each bus $i \in \mathcal{N}_\mathcal{P}$. Buses with $\widetilde{p}_{i}^{\mathcal{P}}>0$ divide loop $\mathcal{P}$ as sub-paths $\{\mathcal{P}_{1},..., \mathcal{P}_{n_r}\}$;\\
	\For{$m=1$ \KwTo $n_r$}
	{
	$\hat{e}_m\leftarrow\text{argmin}_{e\in\mathcal{P}_{m}}\mathbb{E}_{\boldsymbol{\pi}}\left[|p_{e}(\boldsymbol{x}_{\mathcal{E}})|\right]$;
	}
	\For{$e\in \cup_{m=1}^{n_r} \mathcal{K}(\hat{e}_m)$}
	{Solve $\textbf{SOPF}\left(\boldsymbol{\alpha}_{\mathcal{E}\backslash\{e\}}\right)$ to obtain optimal solution $\boldsymbol{x}_{\mathcal{E}\backslash\{e\}}$ and minimum objective value $\widetilde{C}_l(e)$;\\
    $\widetilde{I}_{v}(e)\leftarrow\mathbb{E}_{\boldsymbol{\pi}}\left[F_{vs}\left(\boldsymbol{G}(\boldsymbol{\alpha}_{\mathcal{E}\backslash\{e\}}, \boldsymbol{x}_{\mathcal{E}\backslash\{e\}})\right)\right]$;}
	$e^{*}\leftarrow \text{argmin}_{e\in \cup_{m=1}^{n_r} \mathcal{K}(\hat{e}_m)}\left(k_l\frac{\widetilde{C}_l(e)}{C^{\text{max}}_{l}}+k_{v}\frac{\widetilde{I}_{v}(e)}{I^{\text{max}}_{v}}\right)$;
	\\
	\Return $\boldsymbol{\alpha}^{*}=\boldsymbol{\alpha}_{\mathcal{E}\backslash \{e^{*}\}}$.
\end{algorithm}

We make the following remarks to facilitate understanding of the descriptions in the Algorithm \ref{alg:OSSBR} box.
\begin{itemize}
\item (Line 1) The algorithm is initialized from all the branches being closed, including a redundant branch to open.
\item (Line 3) $\mathcal{P}$ denotes the set of branches forming the single loop in the network, and $\mathcal{N}_{\mathcal{P}}$ is the set of buses in that loop. From any bus $i\in\mathcal{N}_\mathcal{P}$, the expected active power injection into loop $\mathcal{P}$ is
\begin{eqnarray}
	\widetilde{p}_{i}^{\mathcal{P}} := \mathbb{E}_{\boldsymbol{\pi}}\left[ \sum_{\substack{j:i\sim j~{\textnormal{and}}\\j \in \mathcal{N}_{\mathcal{P}}}} p_{ij} \right] = \sum_{w\in\mathcal{W}} \sum_{\substack{j:i\sim j~{\textnormal{and}}\\j \in \mathcal{N}_{\mathcal{P}}}} p_{ij}^w.
	\label{eq:pi_equivalent}
\end{eqnarray}
There are $n_r\geq 1$ buses $i\in\mathcal{N}_{\mathcal{P}}$ satisfying $\widetilde{p}_{i}^{\mathcal{P}}>0$ since the total active power loss in $\mathcal{P}$ must be supplied by at least one source. These $n_r$ buses divide loop $\mathcal{P}$ into $n_r$ sub-paths $(\mathcal{P}_{m},~m=1,...,n_r)$. 

\item (Lines 4--6) $p_e(\boldsymbol{x}_{\mathcal{E}})$ denotes the active power flow on branch $e$ at the optimal solution $\boldsymbol{x}_{\mathcal{E}}$ of $\textbf{SOPF}(\boldsymbol{\alpha}_{\mathcal{E}})$.    
In each sub-path $\mathcal{P}_{m}$,~$m=1,...,n_r$, Line 5 finds the branch $\hat e_m$ that carries the minimum expected absolute value of active power flow.

\item (Line 7) For an arbitrary branch $e=ij\in\mathcal{P}$, denote its upstream branch (in loop $\mathcal{P}$) incident to node $i$ as $e^{\text{up}}$, and downstream branch incident to node $j$ as $e^{\text{down}}$. Based on the expected active branch flow
$\widetilde{p}_{ij}:=\mathbb{E}_{\boldsymbol{\pi}}\left[p_{ij}\right]$ at the optimal solution $\boldsymbol{x}_{\mathcal{E}}$,
we identify a set of candidate branches
\begin{eqnarray}
	&\mathcal{K}(e)&:=
\begin{cases}
	\{e,e^{\text{down}}\}, &\text{if}~\widetilde{p}_{ij}>0~\text{and}~e^{\text{down}}~\text{exists};\\
	\{e,e^{\text{up}}\}, &\text{if}~\widetilde{p}_{ij}<0~\text{and}~e^{\text{up}}~\text{exists};\\
	\{e\}, &\text{otherwise}.
\end{cases}\nonumber
\end{eqnarray}
Starting from Line 7, the algorithm searches in the union of $\mathcal{K}(\hat{e}_{m})$ over all the minimum-flow branches $(\hat e_m,~m=1,...,n_r)$ found in Lines 4--6. 

\item (Line 8) Under switch status $\boldsymbol\alpha_{\mathcal{E}\backslash\{e\}}$ that opens a candidate branch $e$ only and closes all other branches, the algorithm first solves a stochastic OPF problem to minimize the expected total active power loss $\widetilde{C}_l(e)$ only, without considering voltage stability yet (because it lacks the explicit expression of the voltage stability index).


\item \textbf{(Lines 9 and 11)} We highlight them as the key difference from \cite[Algorithm 1]{huang2022improved}. First, Line 9 predicts the expected voltage stability index $\widetilde{I}_{v}(e)$ at the solution obtained in Line 8, using the CNN model from Section \ref{subsec:CNN}. 
Then Line 11 picks a candidate branch $e^*$ that (if opened) can minimize the weighted sum of the expected power loss $\widetilde{C}_l(e)$ and voltage stability index $\widetilde{I}_v(e)$.  
\end{itemize}



\subsection{Deep-Learning-Aided Two-Stage SBR Algorithm}

We now extend the one-stage SBR, Algorithm \ref{alg:OSSBR}, to the two-stage Algorithm \ref{alg:TSSBR} for a general network with $L>1$ redundant branches. The network has $L$ chordless loops $\mathcal{P}^{l}$, $l=1,...,L$. 

\begin{algorithm}[t]
	\caption{Deep-Learning-Aided Two-Stage SBR}\label{alg:TSSBR}
	\SetKwData{Left}{left}\SetKwData{This}{this}\SetKwData{Up}{up}\SetKwFunction{Union}{Union}\SetKwFunction{FindCompress}{FindCompress}\SetKwInOut{Input}{Input}\SetKwInOut{Output}{Output}
	\textbf{First stage:}\\
	Initialize switch status as $\boldsymbol{\alpha}_{\mathcal{E}}:=(\alpha_{ij}=1,\forall ij\in\mathcal{E})$;\\
	Solve $\textbf{SOPF}(\boldsymbol{\alpha}_{\mathcal{E}})$ to get optimal $\boldsymbol{x}_\mathcal{E}=(\boldsymbol{x}^w_\mathcal{E},\forall w\in\mathcal{W})$;\\
	\For{$l=1$ \KwTo $L$}
	{
	$e_{l}^{o}\leftarrow\text{argmin}_{e\in\mathcal{P}^{l}}\mathbb{E}_{\boldsymbol{\pi}}\left[|p_{e}(\boldsymbol{x}_{\mathcal{E}})|\right]$;\\
	Open branch $e_{l}^{o}$;\\
	\If{$e_{l}^{o}$ is a common branch of $\mathcal{P}^{l}$ and $\mathcal{P}^{k}$}{
       Update loop $\mathcal{P}^{k}$ with $e_{l}^{o}$ open;
      }	
	}
    $\mathcal{E}^o\leftarrow\{e_{l}^{o},~l=1,...,L\}$;\\
	\textbf{Second stage:} Preset $N = N_{\text{max}}$;\\
	\For{$n=1$ \KwTo $N_{\textnormal{max}}$}{
	\For{$l=1$ \KwTo $L$}
	{
	Define a set of open branches $\mathcal{E}^o_{n,l}:=\mathcal{E}^o\backslash\{e_l^o\}$;\\ 
	Call \textbf{Algorithm \ref{alg:OSSBR}} with initial switch status $\boldsymbol{\alpha}_{\mathcal{E}\backslash\mathcal{E}^o_{n,l}}$ to obtain an optimal open branch $e_{n,l}^*$ and the corresponding minimum objective value $\Gamma_{n,l}^*$ (in Line 11, Algorithm \ref{alg:OSSBR});\\
	$e_l^o\leftarrow e_{n,l}^*$, which also updates $\mathcal{E}^o$;
	}
	\If{$\Gamma_{n,1}^*=\Gamma_{n,2}^*=...=\Gamma_{n,L}^*$}{
	    $N\leftarrow n$; break;
      }
	}
	$(n^*,l^*)\leftarrow\text{argmin}_{\substack{n=1,...,N\\l=1,...,L}}~ \Gamma_{n,l}^*$;\\
	\Return $\boldsymbol{\alpha}^{*}=\boldsymbol{\alpha}_{\mathcal{E}\backslash\left(\mathcal{E}^o_{n^*,l^*} \cup \{e_{n^*,l^*}^*\}\right)}$.
\end{algorithm}

The first stage (Lines 1--11) of Algorithm \ref{alg:TSSBR} is the same as \cite[Algorithm 2]{huang2022improved}. It starts from all the branches being closed, solves a stochastic OPF problem, and sequentially opens a branch carrying the minimum expected active power flow in each of the $L$ loops.  
If such an open branch simultaneously lies in two loops, then opening it will also require updating the other loop. 
The $L$ open branches, collected in a set $\mathcal{E}^o$, serve as the initial candidate branches to open for the next stage. 



The second stage of Algorithm \ref{alg:TSSBR} inherits the iterative close-and-open idea from \cite[Algorithm 2]{huang2022improved}: in each of the $L$ (inner) iterations, it opens all but one branches in the set $\mathcal{E}^o$ (Line 15). It then calls Algorithm \ref{alg:OSSBR} to find the single redundant branch that remains to open (Line 16), and uses that branch to update the set $\mathcal{E}^o$ (Line 17). 
Besides the CNN-based prediction of voltage stability already in Algorithm \ref{alg:OSSBR}, another difference from \cite{huang2022improved} lies in the outer iteration with index $n$ (Line 13). 
After traversing (and perhaps updating) all the $L$ branches in $\mathcal{E}^o$, the next outer iteration will repeat the close-and-open procedure on $\mathcal{E}^o$, until the minimum objective value of \textbf{SDNR-VS} gets steady as the close-and-open procedure continues (Line 19) or the limit $N_{\text{max}}$ of outer iterations is reached. 
Algorithm \ref{alg:TSSBR} ultimately returns the best radial topology found throughout all the outer and inner iterations.

The case studies in Section \ref{Sec:CaseStudy} will show that Algorithm \ref{alg:TSSBR} has satisfactory computational efficiency.
Indeed, the CNN-based prediction of voltage stability takes moderate extra time than \cite[Algorithm 2]{huang2022improved}, and a high-quality \textbf{SDNR-VS} solution is usually obtained within $N = 2$ outer iterations. 


\section{Case Studies} \label{Sec:CaseStudy}

\subsection{Experimental Setup}

To validate the proposed method, we conduct numerical experiments on the IEEE 33-bus and 123-bus distribution network models that were also used in \cite{huang2022improved}. 
Each network has one substation bus. The 33-bus network has four renewable generation buses (each connected to a small wind turbine and a solar panel) and five redundant branches; the 123-bus network has six renewable generation buses and four redundant branches (after removing the branch between bus 17 and bus 52). 
Other non-substation buses connect to a load each. 
We use the load, wind and solar generation data in Germany from January 2018 to June 2020 in hourly resolution \cite{OPSDP2022} and scale them to fit the test network capacities. After getting the active power of loads and renewable generations, we obtain their reactive power using fixed power factors \cite{zhan2020switch}. 
The scaled data from January 2018 to March 2020 are sampled for training and testing the CNN-based voltage stability prediction model, and the remaining data from April to June 2020 are clustered into a certain number of scenarios (depending on the test case) using k-medoids \cite{huang2021distribution} to validate the proposed method.

We split the data into 70\% for training and 30\% for testing the CNN model. 
Each data sample consists of an input $\boldsymbol{G}(\boldsymbol{\alpha},\boldsymbol{x})$ (network state) and an output $I_{v}$ (the smallest singular value $\delta_{\text{min}}$ of the power-flow Jacobian for steady-state voltage stability; or the RVSI for short-term voltage stability): 

\begin{itemize}
\item For steady-state voltage stability, we solve the OPF problems given all the possible radial configurations (there are 33,913 and 42,658 of them, respectively, for the 33-bus and 123-bus networks) and the sampled loads and renewable generations (discarding those without feasible power flow solutions). 
The stability index $\delta_{\text{min}}$ is  calculated for each of the OPF solutions using singular value decomposition. 

\item For short-term voltage stability, we connect the distribution substation to the IEEE 39-bus transmission network, at transmission bus 37 as the default point of common coupling (PCC). 
Among all the possible radial configurations of the distribution network, 20,000 cases are sampled. For each sample, time-series simulation is conducted in a 5-second window for a three-phase short-circuit fault at the PCC. 
An RVSI is calculated from the time-series voltage trajectories in each simulation. 
The induction motor models for the simulation are generated from the field measurements in \cite{renmu2006composite}, and the ZIP loads are uniformly randomly sampled within [1, 2] per unit. In case a voltage collapse occurs in a simulation, the corresponding RVSI is set to a very large value.
\end{itemize}

Our CNN model is composed of four convolutional layers with 8, 16, 32, and 64 filters, respectively; each filter is $3\times 3$; the learning rate is $1\times 10^{-3}$; the maximum epoch is 30; the mini-batch size is 20; and the dropout rate is 0.2. 
The experiments are run on a 64-bit MacBook with 
8-core CPU and 32GB RAM. 
We use the Matpower Interior Point Solver for OPF, the PSAT for time-series simulation, and the MATLAB Deep Learning Toolbox to build the CNN.

The following methods are compared in the case studies: 
\begin{itemize}
\item \textbf{The proposed method:} the deep-learning-aided two-stage SBR for \textbf{SDNR-VS}, i.e., Algorithm \ref{alg:TSSBR}. The weighting factors in objective \eqref{eq:SDNR_obj} are set as $k_{l} = 0.5$ and $k_{v}=\pm 0.5$ (negative when $I_v^w = \delta_{\text{min}}$ for steady-state stability and positive when $I_v^w = \text{RVSI}$ for short-term stability).

\item \textbf{Method 1:} using the mixed-integer nonlinear program solver Gurobi to solve SDNR \textit{without considering voltage stability}. We use Method 1 as a benchmark by defining the relative errors of other methods compared to it: $\eta_{C_l}$ for the expected total active power loss, $\eta_{\delta_{\text{min}}}$ for the expected $\delta_{\text{min}}$, and $\eta_{\text{RVSI}}$ for the expected RVSI. 
Note that a higher $\eta_{\delta_{\text{min}}}$ indicates better steady-state voltage stability and a lower $\eta_{\text{RVSI}}$ indicates better short-term voltage stability.
	
\item \textbf{Method 2:} the two-stage SBR for SDNR \textit{without considering voltage stability}, i.e., \cite[Algorithm 2]{huang2022improved}. 

\item \textbf{Method 3:} replacing the CNN model in the proposed method with singular value decomposition for steady-state voltage stability and simulation-based evaluation for short-term voltage stability.
It thus provides accurate voltage stability evaluation at high computational expenses. 
\end{itemize}

\subsection{Accuracy of CNN Prediction for Voltage Stability}

For the $i$-th data sample, let $I_{v,i}$ and $\hat{I}_{v,i}$ denote the real and CNN-predicted voltage stability indices, respectively. 
It is actually the order, in which the voltage stability indices of different configurations are ranked, that matters most.    
Therefore, we measure the accuracy of the CNN prediction by the \textit{consistency} between the real and predicted orders \cite{huang2021distribution}:
\begin{eqnarray}
\text{Consistency}:=\frac{2\sum_{i=1}^{H-1}\sum_{j=i+1}^{H}c(i,j)}{H(H-1)}\times 100\% \nonumber
\end{eqnarray}
where for each pair $i,j=1,...,H$ of data samples, $c(i,j)=1$ if the order between $(I_{v,i},I_{v,j})$ is the same as that between $(\hat{I}_{v,i},\hat{I}_{v,j})$, and $c(i,j)=0$ otherwise.

\begin{table}
	\centering
	\caption{Consistency of CNN Prediction for Voltage Stability (10 Runs)}
\renewcommand\arraystretch{1.2}	
	\begin{tabular}{c|c|ccc}		
		\hline
		\hline 
		Network 
		&Voltage stability index&mean (\%) &max (\%) &min (\%)\\
		\hline 
		\multirow{2}{*}{33-bus}
		&$\delta_{\text{min}}$ (steady-state)  &94.9 &95.1 &94.6\\
		&RVSI (short-term)&94.6 &94.7 &94.4\\
		\hline 
		\multirow{2}{*}{123-bus}
		&$\delta_{\text{min}}$ (steady-state) &97.7 &97.9 &97.6\\
		&RVSI (short-term)&98.8 &98.8 &98.7\\
		\hline
		\hline  	
	\end{tabular}  \label{tab:CNNReliability}
\end{table}

\begin{figure}
	\centering	\includegraphics[width=0.9\columnwidth]{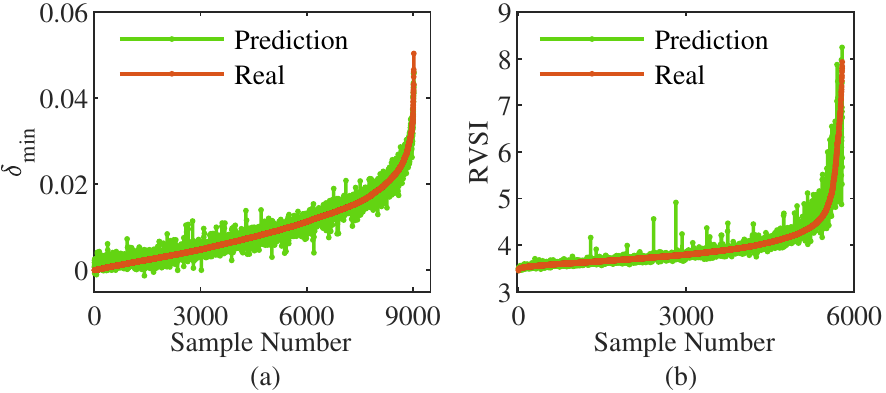}
	\caption{Comparison of real and CNN-predicted (a) $\delta_{\text{min}}$ for steady-state voltage stability and (b) RVSI for short-term voltage stability. Samples are taken from the 33-bus network and ranked in an ascending order of the voltage stability indices.}
\label{fig:CNNReliability}
\end{figure}

We run ten CNN training processes with different random samples, which all result in good consistency as Table~\ref{tab:CNNReliability} shows. 
Figure~\ref{fig:CNNReliability} compares the real and predicted voltage stability indices for samples from the 33-bus network. It shows good accuracy of the proposed CNN prediction model.

\subsection{Voltage Stability Enhancement}

\begin{table*}
	\centering
	\caption{Power Loss and Voltage Stability: $|\mathcal{W}|=5$, Different Networks, Renewable Penetration Levels, and Methods}
	\renewcommand\arraystretch{1.2}	
	\begin{tabular}{c|c|c|ccc|ccc|ccc|ccc}	
		\hline
		\hline 
		\multirow{3}{*}{Network}&\multirow{3}{*}{$k_r$} 
		&\multirow{3}{*}{Method}  
		&\multicolumn{6}{c|}{Method 3 and the proposed: enhance $\delta_{\text{min}}$}
		&\multicolumn{6}{c}{Method 3 and the proposed: enhance RVSI}\\
		\cline{4-15}
		& &
		&\multicolumn{3}{c|}{\makecell{$\eta_{C_l}$ (\%)}}  &\multicolumn{3}{c|}{\makecell{$\eta_{\delta_{\text{min}}}$ (\%)}}
		&\multicolumn{3}{c|}{\makecell{$\eta_{C_l}$ (\%)}}  &\multicolumn{3}{c}{\makecell{$\eta_{\text{RVSI}}$ (\%)}}\\
		\cline{4-15}
		&& &mean  &max &min &mean  &max &min  &mean  &max &min &mean  &max &min \\
		\hline 
		\multirow{9}{*}{33-bus}&\multirow{3}{*}{1.0}
		&Method 2 &-1.07 &0.21 &-5.82 &36.06 &128.35 &-14.30
		&-1.07 &0.21 &-5.82 &0.23 &28.29 &-40.27
		\\
		&&Method 3 &-1.03 &0.21 &-5.82 &42.66 &128.35 &-14.30
		&-0.50 &0.64 &-5.37 &-6.69 &0 &-41.79
		\\
		&&\textbf{Proposed} &\textbf{-1.03} &\textbf{0.21} &\textbf{-5.82} &\textbf{42.66} &\textbf{128.35} &\textbf{-14.30}
	    &\textbf{-0.31} &\textbf{2.24} &\textbf{-5.82} &\textbf{-6.09} &\textbf{2.28} &\textbf{-41.81}
		\\
		\cline{2-15} 
		&\multirow{3}{*}{1.5}
		&Method 2 &-0.86 &0.22 &-4.51 &51.15 &127.96 &-0.28
		&-0.86 &0.22 &-4.51 &-5.78 &23.12 &-39.42
		\\
		&&Method 3 &-0.71 &1.81 &-4.51 &51.97 &127.96 &-0.28
		&-0.31 &1.39 &-3.69 &-7.92 &11.91 &-39.42
		\\
		&&\textbf{Proposed} &\textbf{-0.86} &\textbf{0.22} &\textbf{-4.51} &\textbf{51.15} &\textbf{127.96} &\textbf{-0.28}
		&\textbf{3.77} &\textbf{12.79} &\textbf{-1.20} &\textbf{-9.42} &\textbf{3.11} &\textbf{-41.13}
		\\
		\cline{2-15}
		&\multirow{3}{*}{2.0}
		&Method 2 &-1.29 &0.57 &-5.87 &41.01 &148.93 &-44.37
		&-1.29 &0.57 &-5.87 &-3.43 &36.53 &-40.61
		\\
		&&Method 3 &-1.09 &0.65 &-5.87 &55.19 &210.16 &-15.32
		&-0.16 &4.55 &-4.09 &-10.90 &2.08 &-40.61
		\\
		&&\textbf{Proposed} &\textbf{-1.32} &\textbf{0.65} &\textbf{-5.87} &\textbf{51.93} &\textbf{184.17} &\textbf{-15.32}
		&\textbf{1.33} &\textbf{4.95} &\textbf{-4.09} &\textbf{-8.59} &\textbf{62.72} &\textbf{-41.51}
		\\
		\hline  	
		\multirow{6}{*}{123-bus}&\multirow{2}{*}{1.0}
		&Method 2 
		&-0.03 &0 &-0.16 &1.30 &8.66 &-0.77
		&-0.03 &0 &-0.16 &-2.11 &0.13 &-50.91
		\\
		&&\textbf{Proposed} 
		&\textbf{-0.02} &\textbf{0.03} &\textbf{-0.16} &\textbf{1.38} &\textbf{8.66} &\textbf{-0.77}
		&\textbf{0.00} &\textbf{0.03} &\textbf{-0.14} &\textbf{-2.21} &\textbf{0.00} &\textbf{-50.97}
		\\
		\cline{2-15} 
		&\multirow{2}{*}{1.5}
		&Method 2 
		&-0.03 &0.11 &-0.22 &2.13  &10.63 &-0.22
		&-0.03 &0.11 &-0.22 &-4.30 &0.11 &-51.42
		\\
		&&\textbf{Proposed} 
		&\textbf{-0.03} &\textbf{0.12} &\textbf{-0.21} &\textbf{2.21}  &\textbf{10.44}  &\textbf{-0.22}
		&\textbf{-0.02} &\textbf{0.11} &\textbf{-0.22} &\textbf{-4.38}  &\textbf{-0.03}  &\textbf{-51.43}
		\\
		\cline{2-15} 
		&\multirow{2}{*}{2.0}
		&Method 2 
		&0.00 &0.27 &-0.63 &2.41  &9.43 &-1.99
		&0.00 &0.27 &-0.63 &-8.59 &0.17 &-52.70
		\\
		&&\textbf{Proposed} 
		&\textbf{0.00} &\textbf{0.27} &\textbf{-0.63} &\textbf{2.55}  &\textbf{9.43}  &\textbf{-0.21}
		&\textbf{-0.01} &\textbf{0.20} &\textbf{-0.63} &\textbf{-8.61}  &\textbf{0.07}  &\textbf{-52.70} 
		\\
		\hline  	
		\hline  	
	\end{tabular}
 \label{tab:B33-123-vs}
\end{table*}


Table~\ref{tab:B33-123-vs} shows
the results of different methods applied to the 33-bus and 123-bus networks (except that Method 3 is not applied to the 123-bus network due to its huge computation time). 
Different renewable penetration levels are realized by scaling the renewable generation capacities with a factor $k_r$. 
The numbers are the relative errors of different methods in power loss and voltage stability indices compared to Method 1. 
Specifically, in the first six columns $\delta_{\text{min}}$ is enhanced, and in the last six columns the RVSI is enhanced, by the proposed method and Method 3 for \textbf{SDNR-VS}. 
Although these voltage stability indices are excluded from the objectives of Methods 1 and 2, they can still be calculated after the SDNR problems are solved.
The mean, maximum, and minimum refer to those of the 24 records in a day, each for an hour.

\begin{figure}
	\centering
\includegraphics[width=0.9\columnwidth]{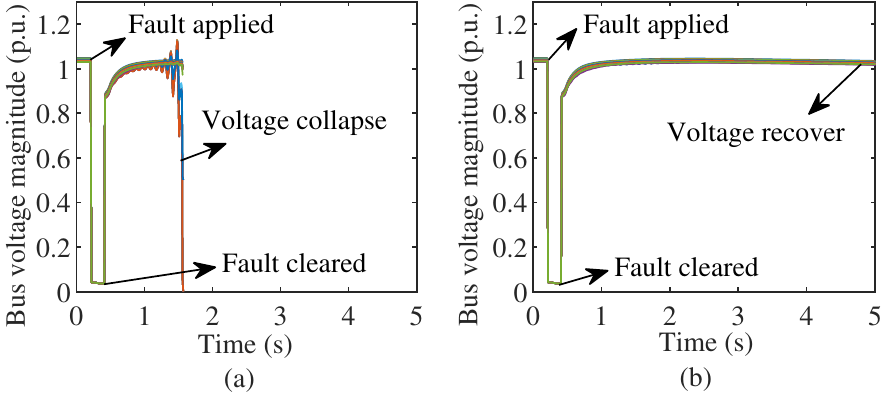}
	\caption{Voltages after a fault at different buses (in different colors) in the 33-bus network, under (a) Method 2 without enhancing voltage stability and (b) the proposed method. Renewable penetration level $k_r = 2$.}	\label{fig:B33_STVS_kr2}
\end{figure}

The main observations from Table \ref{tab:B33-123-vs} are:
\begin{itemize}
\item In most cases, the proposed method and Method 3 improve the voltage stability index $\delta_{\text{min}}$ or RVSI, compared to Methods 1 and 2 (especially Method 1) that do not include voltage stability in their objectives. 
Indeed, the SDNR decisions from Method 2 may lead to voltage collapses, which can be prevented by the proposed method. Such an example is shown in Figure~\ref{fig:B33_STVS_kr2}: the 5-second record of bus voltages after a fault, in the 33-bus network.  

\item The improvement in voltage stability is mostly achieved with minor difference in power loss. The worst (highest) mean increase of power loss is $\eta_{C_l}=3.77\%$ when renewable penetration $k_r=1.5$ and the proposed method enhances RVSI. This is still an acceptable trade-off.   

\item The proposed method and Method 3 show similar results in most cases in the 33-bus network. This further verifies the accuracy of the CNN-based prediction of voltage stability. 
There is a rare case where the proposed method gets much worse RVSI (maximum $\eta_{\text{RVSI}}=62.72\%$) than Method 3, which is likely due to the large prediction error of the CNN under a high renewable penetration level $k_r=2$. Note that our CNN model is trained with data sampled under $k_r=1$ only. Hence it is reasonably conjectured that training different CNNs under different renewable penetration levels may further improve the accuracy of prediction and the performance of the proposed method.      
\end{itemize}

\begin{figure*}
\centering
\subfigure[Total active power loss (enhancing $\delta_{\text{min}}$)]{
\includegraphics[width=0.9\columnwidth]{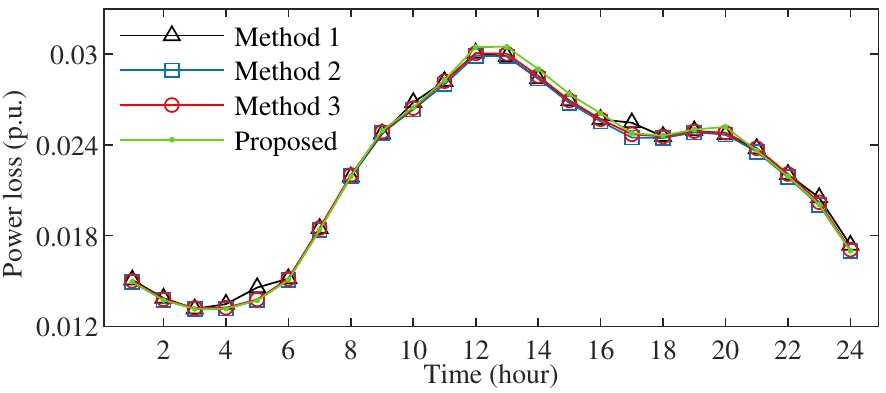}
}\label{fig:B33_svs_pls}
\subfigure[Stead-state voltage stability index ($\delta_{\text{min}}$)]{
\includegraphics[width=0.9\columnwidth]{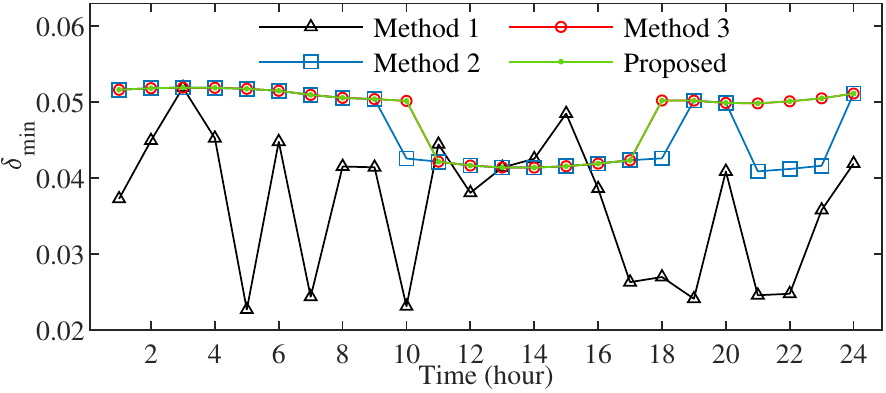}
}\label{fig:B33_svs_idx}
\subfigure[Total active power loss (enhancing RVSI)]{
\includegraphics[width=0.9\columnwidth]{SDNRSTVS_kres1_ploss-eps-converted-to.pdf}
}\label{fig:B33_stvs_pls}
\subfigure[Short-term voltage stability index (RVSI)]{
\includegraphics[width=0.9\columnwidth]{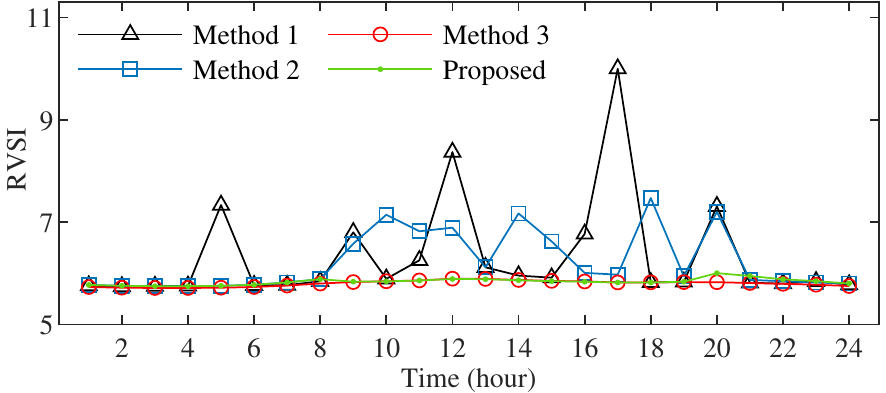}
}\label{fig:B33_stvs_idx}
\caption{One-day results of different methods in the 33-bus network under renewable penetration level $k_{r} = 1$. The proposed method and Method 3 solve \textbf{SDNR-VS} to enhance steady-state voltage stability (a), (b) or short-term voltage stability (c), (d).} 
\label{fig:B33_vs}
\end{figure*}


The observations above are supplemented by Figure \ref{fig:B33_vs}, which displays the 24-hour power loss and voltage stability indices when different methods are applied. It verifies that the proposed method can significantly improve voltage stability with minor (if any) increase in power loss, and its CNN-based performance is very close to Method 3 that makes ``ground-truth'' voltage stability evaluation.


\begin{table*}
	\centering
	\caption{Power Loss and Voltage Stability: 33-Bus Network, $k_r = 1.0$, Different $|\mathcal{W}|$, Different Methods}
	\renewcommand\arraystretch{1.2}	
	\begin{tabular}{c|c|ccc|ccc|ccc|ccc}	
		\hline
		\hline 
		\multirow{4}{*}{$|\mathcal{W}|$} 
		&\multirow{3}{*}{Method}  
		&\multicolumn{6}{c|}{The proposed method: enhance $\delta_{\text{min}}$}
		&\multicolumn{6}{c}{The proposed method: enhance RVSI}\\
		\cline{3-14}
		& 
		&\multicolumn{3}{c|}{\makecell{$\eta_{C_l}$ (\%)}}  &\multicolumn{3}{c|}{\makecell{$\eta_{\delta_{\text{min}}}$ (\%)}}
		&\multicolumn{3}{c|}{\makecell{$\eta_{C_l}$ (\%)}}  &\multicolumn{3}{c}{\makecell{$\eta_{\text{RVSI}}$ (\%)}}
		\\
		\cline{3-14}
		& &mean  &max &min &mean  &max &min  &mean  &max &min &mean  &max &min \\
		\hline 
		\multirow{2}{*}{20}
		&Method 2 
		&-2.44 &0.33 &-6.64 &60.71 &174.93 &1.63
		&-2.44 &0.33 &-6.64 &5.60 &238.90 &-36.07
		\\
		&\textbf{Proposed} 
		&\textbf{-2.40} &\textbf{0.33} &\textbf{-6.64} &\textbf{65.67} &\textbf{174.93} &\textbf{15.30}
		&\textbf{-1.78} &\textbf{1.02} &\textbf{-6.47} &\textbf{-12.02} &\textbf{15.69} &\textbf{-37.35}
		\\
		\hline 
		\multirow{2}{*}{40}
		&Method 2 
		&-14.54 &-2.76 &-55.31 &185.05  &1056.10 &44.83
		& -14.54 &-2.76 &-55.31 & -31.78 &23.13 &-75.45
		\\
		&\textbf{Proposed} 
		&\textbf{-14.51} &\textbf{-2.76} &\textbf{-55.31} &\textbf{192.89}  &\textbf{1056.10}  &\textbf{76.02}
		&\textbf{-13.90} &\textbf{-0.68} &\textbf{-55.31} &\textbf{-37.21}  &\textbf{3.17}  &\textbf{-77.27} 
		\\
		\hline  	
		\hline  	
	\end{tabular} \label{tab:B33_svs_scenario}
\end{table*}

We perform extended experiments on the 33-bus network with increasing numbers of uncertainty scenarios $|\mathcal{W}|$. The results are shown in Table~\ref{tab:B33_svs_scenario}. As $\mathcal{|W|}$ increases, the proposed method obtains better solutions with smaller power loss and better voltage stability. 
Compared to Method 2, the proposed method makes significant improvement in voltage stability index $\delta_{\text{min}}$ or RVSI, with small changes in power loss.

\begin{figure}
	\centering
\includegraphics[width=0.9\columnwidth]{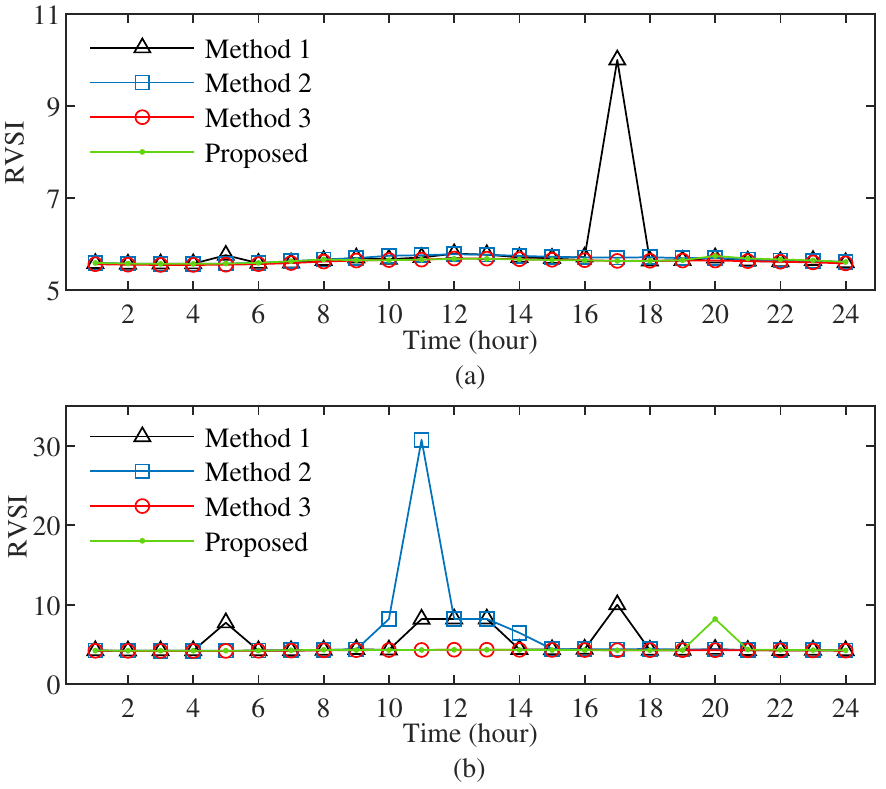}
	\caption{RVSI of 33-bus networks reconfigured by different methods and connected to (a) bus 27 or (b) bus 20 of the transmission network.}
	\label{fig:B33_PCCs}
\end{figure}

Figure~\ref{fig:B33_PCCs} compares the RVSI of the 33-bus networks reconfigured by different methods and connected to different PCC buses in the transmission. 
An RVSI is calculated for a three-phase short-circuit fault in each of the 24 hours. 
Both Methods 1 and 2 (without enhancing voltage stability) encounter voltage collapses (indicated by $\text{RVSI}\geq 10$) at some time, which are fixed by the proposed method. Still, the proposed method achieves similar performance to Method 3 that makes accurate RVSI evaluation from the time-series voltage record.

\begin{table}
	\centering
	\caption{Computation Time of the Proposed Method (in seconds)}
	\renewcommand\arraystretch{1.2}	
	\begin{tabular}{c|c|ccc}	
		\hline
		\hline 
		Network	&$|\mathcal{W}|$ &Method 1 &Enhance $\delta_{\text{min}}$ &Enhance RVSI  \\
		\hline 
        \multirow{3}{*}{33-bus} 
        &5 &166.2 &4.1 &5.1 \\
		&20  &316.9 &15.9 &19.7 \\	
		&40 &901.3 &31.8 &40.4 \\
		\hline
		\multirow{3}{*}{123-bus} 
		&5 &494.1 &5.1 &4.9\\
		&20  &832.4 &24.2 &20.0\\	
		&40 &3333.9 &56.2 &46.6\\
		\hline
		\hline  	
	\end{tabular} \label{tab:B33_time}
\end{table}

\subsection{Computational Efficiency}

Table~\ref{tab:B33_time} compares the average computation time of Method 1 (using Gurobi) and the proposed method over 24 hours, under different uncertainty  scenario numbers $|\mathcal{W}|$. Compared to Method 1, the proposed method speeds up the computation by at least an order of magnitude. 

\section{Conclusion}\label{sec:conclusion}
We proposed a deep learning method to solve stochastic distribution network reconfiguration (SDNR) with voltage stability enhancement. 
A convolutional neural network (CNN) model for voltage stability prediction is integrated into successive branch reduction (SBR) algorithms to search for a radial topology with optimized performance in terms of power loss reduction and voltage stability enhancement. 
Numerical experiments on the IEEE 33-bus and 123-bus network models verified that the proposed method can significantly improve steady-state or short-term voltage stability with a minor compromise in power loss optimization. Moreover, it speeds up the solution process of SDNR by at least an order of magnitude, compared to a classic mixed-integer nonlinear program solver.   

In the future, we plan to formally analyze the optimality of the proposed SBR algorithms (which are currently heuristics), incorporating the CNN-induced prediction errors in voltage stability indices. 
Another objective of our ongoing work is to co-optimize the placement of volt/var control resources (capacitor banks, smart inverters, D-STATCOMs, etc.) in synergy with the topology reconfiguration.
\bibliographystyle{IEEEtran}
\bibliography{references}

\end{document}